\documentclass[11pt]{article}
\usepackage{amsmath}
\usepackage{amscd}
\usepackage{epsfig}

\usepackage[latin1]{inputenc}

\input{amssym}

\newcommand{\color}[6]{}
\newcommand\QQ{\hbox{I\kern-.53em\hbox{Q}}}
\newcommand\qed{\hfill$\sqcap\kern-8.0pt\hbox{$\sqcup$}$}
\newcommand\NN{\hbox{I\kern-.2em\hbox{N}}}
\newcommand\RR{\hbox{I\kern-.2em\hbox{R}}}
\newcommand\sRR{{\sl \hbox{I\kern-.2em\hbox{R}}}}
\newcommand{\PP}{{\bf P}^k}

\newcommand{\pp}{{\bf P}^1}

\newcommand\ZZ{{{\rm Z}\kern-.28em{\rm Z}}}
\newcommand\proof{\noindent{\em{Proof}.\ }}

\newcommand{\Bif}{{\cal B}{\textrm{\scriptsize if}}}
\newcommand{\Bc}{T_{\textrm{\scriptsize bif}}}
\newcommand{\Bm}{\mu_{\textrm{\scriptsize bif }}}

\newcommand{\car}{\heartsuit}
\newcommand{\Car}{\bf{\heartsuit}}

\newtheorem{theo}{Theorem}[section]
\newtheorem{prop}[theo]{Proposition}
\newtheorem{lem}[theo]{Lemma}

\newtheorem{cor}[theo]{Corollary}

\newtheorem{defi}[theo]{Definition}
\newtheorem{defiprop}[theo]{Definition-Proposition}

\newcommand{\la}{\lambda}

\numberwithin{equation}{section}

\begin{document}

\date{}
\title{Lyapunov exponents, bifurcation currents and laminations in bifurcation loci}

\vskip0.5cm

\author{Giovanni Bassanelli and Fran\c{c}ois Berteloot
}

\vskip0.5cm

\maketitle

\begin{abstract}
Bifurcation loci in the moduli space of degree $d$ rational maps are shaped by the hypersurfaces
defined by the existence of a cycle of period $n$ and multiplier $0$ or $e^{i\theta}$.
Using potential-theoretic arguments, we establish two equidistribution properties for these hypersurfaces
with respect to the bifurcation current. To this purpose we first establish approximation formulas for the Lyapunov
function.
In degree $d=2$, this allows us to build holomorphic motions and show that the bifurcation locus
has a lamination structure in the regions where an attracting basin of fixed period exists.
\end{abstract}

\vskip0.5cm

{\footnotesize Giovanni Bassanelli, Dipartimento di Matematica, 
Universit\`a di Parma, Parco Area delle Scienze,
Viale Usberti 53/A I--43100 Parma, Italia.
{\em Email: giovanni.bassanelli@unipr.it}

Fran\c cois Berteloot, Universit\' e Paul Sabatier MIG. 
Institut de Math\'ematiques de Toulouse. 
31062 Toulouse
Cedex 9, France.
{\em Email: berteloo@picard.ups-tlse.fr}}
\medskip

\noindent Math Subject Class: 37F45; 37F10

\section{Introduction}
A bifurcation is usually said to occur when the dynamical behaviour of a map changes drastically under perturbation. 
For rational maps on the 
Riemann sphere, the study of bifurcations started with the
seminal paper of Ma\~n\'e, Sad and Sullivan \cite{MSS} and, since, constitutes a large field of research. Let us recall that, in particular, this work 
characterizes the bifurcation locus as the subset of the parameter space where the Julia set does not move continuously.\\
In any parameter space $M$ of some holomorphic family $\left(f_{\la}\right)_{\la\in M}$ or in the moduli space ${\cal M}_d$ of degree $d$ rational maps, 
the subsets $Per_n(w)$ of rational maps possessing a cycle of exact period $n$ and multiplier $w\in {\bf C}\setminus \{1\}$ turn out to be hypersurfaces. 
As the  
Ma\~n\'e-Sad-Sullivan theory shows, the bifurcation locus is shaped by the hypersurfaces $Per_n(0)$ and $Per_n(e^{i\theta})$; it actually coincides with the 
closure of the union over $n$ and $\theta$ of the hypersurfaces  $Per_n(e^{i\theta})$. Moreover, the work of Milnor \cite{Mi} reveals that
the global structure of the bifucation locus is intimately related to
the behaviour of these hypersurfaces in a convenient 
compactification of 
${\cal M}_d$. We aim, in the present paper, to combine the approach of Milnor with 
the point of view of bifurcation currents. 
As we shall now explain, these currents allow to
 both exploit the properties of the Lyapunov function and general measure or potential-theoretic methods.\\

The bifurcation locus of any holomorphic family $\left(f_{\la}\right)_{\la\in M}$, as well as the bifurcation locus of a moduli space ${\cal M}_d$,
 carries a closed, positive $(1,1)$-current. This current, which has been introduced by DeMarco \cite{DeM1}, 
is called the bifurcation current and denoted $\Bc$. Both $\Bc$ and its powers 
$\Bc^k$ ($k\le dim_{\bf C}M$) are extremely convenient tools for the study of measurable or complex 
analytic properties of the bifurcation locus and have been exploited in several recent works (\cite{BB},\cite{Ph},\cite{DF},\cite{Du}).  
A crucial fact is  
that the bifurcation current admits both the Lyapunov function $L(\la)$ and the sum of values of the Green function on the critical points 
 as potentials (see \cite{DeM2} or \cite{BB}):
$$\Bc=dd^c L(\la)=dd^c \sum G_{\la}(c_\la).$$ 
We recall that $L(\la)$ denotes the Lyapunov exponent of $f_{\la}$ with respect to its maximal entropy measure.
The above formula actually enlights the double nature of the bifurcation current: 
$\Bc$ may either detect the instability of repelling cycles
or the instability of critical orbits. This 
reflects the essence of the classical Ma\~n\'e-Sad-Sullivan theory on bifurcations (see
\cite{BB}, Section 5). Our point of view, initiated in \cite{BB}, is based on the identity $\Bc=dd^c L(\la)$.
\\

Let us now describe the content of the paper.
We first investigate how the hypersurfaces $Per_n(0)$ and $Per_n(e^{i\theta})$
accumulate the bifurcation current and establish  
the following equidistribution formulas in a quite general context (see theorem \ref{appT} and, for similar results regarding $\Bc^k$,
theorem \ref{AppMu}):

\begin{eqnarray*}
\Bc&=&\lim_n d^{-n}[Per_n(0)]\\
\Bc&=& \lim_n \frac{d^{-n}}{2\pi}\int_{0}^{2\pi}[Per_n(e^{i\theta})]d\theta.
\end{eqnarray*}

\noindent It is important to stress here that these formulas are deduced from the identity $\Bc=dd^c L(\la)$ and the following property of the Lyapunov exponent
($RP(f_{\la},n)$ is the set of repulsiv points of $f_{\la}$ of exact period $n$):
$$L(\la)=\lim_n d^{-n}\sum_{RP(f_{\la},n)}\log\vert f_{\la}'\vert.$$
A similar fact has been proved for holomorphic endomorphisms of ${\PP}$ (\cite{BDM}). Although the proof is definitely easier when $k=1$,
 the above formula seems to be new even in that case.

In the remaining of the paper we work in the moduli space ${\cal M}_2$ of quadratic rational maps. 
As it has been shown by Milnor \cite{Mi}, we may identify ${\cal M}_2$ with ${\bf C}^2$.  
Our goal is to
use the above formulas to see how the distribution of the hypersurfaces $Per_{n}(0)$ determines the shape of the bifurcation locus.\\
Besides the equidistribution formulas, we introduce two other tools which may be of independant interest. The first one is a parametrization
``\`a la Douady-Hubbard'' of certain hyperbolic components of ${\cal M}_2$ (see theorem \ref{thUni}). This parametrization yields
a holomorphic motion of the hyperbolic components of any curve $Per_{n_0}(0)$. 
The second is an extension principle for special kind of holomorphic motions in ${\bf P}^2$. This principle is comparable with
the classical $\la$-lemma 
for holomorphic motions in the Riemann sphere. It allows us
to obtain the motion of the bifurcation locus 
 as a limit of motions of hyperbolic components (see theorem \ref{thext}).
 
Our main result describes the bifurcation locus $\Bif_{n_0}$ of the open region $U_{n_0}\subset {\cal M}_2$ consisting of ``parameters''
 which do have an attracting
cycle of period $n_0$. We compare $\Bif_{n_0}$ with the bifurcation locus $\Bif_{n_0}^{\;c}$ of the ``central curve'' $Per_{n_0}(0)$ (see theorem \ref{thMouvBif}
for a precise statement). We show that $\Bif_{n_0}$ is a lamination whose transverse measure is the bifurcation measure $\mu_{n_0}^{\;c}$ on the 
central curve $Per_{n_0}(0)$.
We actually construct, in a very natural way, a
 holomorphic motion $\sigma:\Bif_{n_0}^{\;c}\times\Delta \longrightarrow \Bif_{n_0}$ such that:
$$\Bc{\arrowvert_{U_{n_0}}} = \int_{\Bif_{n_0}^{\;c}}\left[\sigma\left(\la,\Delta\right)\right]\;\mu_{n_0}^{\;c}.$$
When $n_0=1$, our result is sharper and also easier to visualize. The open region $U_1$ is the subset of ${\cal M}_2$ which consists of parameters
having an attracting fixed point. The central curve $Per_1(0)$ is a line in ${\bf C}^2={\cal M}_2$. This line is the family
of quadratic polynomials on which the bifurcation locus
 is 
the boundary of the Mandelbrot set ${M}_2$. The bifurcation measure is the harmonic measure of ${M}_2$
(see corollary \ref{mandel}). Our result says, in particular, that the bifurcation locus in $U_1$ is
obtained by moving holomorphically the boundary of the Mandelbrot set. 
Using Slodkowski's theorem, we extend our
holomorphic motion to the full line $Per_1(0)$. This shows that every stable component moves holomorphically. Consequently, we observe
that if some non-hyperbolic stable component
would exist on the line of quadratic polynomials, such a component would also exist within quadratic rational maps (see theorem \ref{farf}).\\

We now end this intoduction with a few words on related works.
In their paper (\cite{DF}, Theorems 1 and 4.2) Dujardin and Favre have studied the distribution of hypersurfaces
of $M$ defined by the pre-periodicity of a critical point. They have shown that the bifurcation
current is the limit of  suitably normalized currents of integration on these hypersurfaces. The first assertion of our
theorem \ref{appT} corresponds to the case of a periodic critical point. It is interesting to note that the approaches 
are based on the two different interpretations of the bifurcation current which we have previously described. Dujardin and Favre
use $\Bc=dd^c \sum G_{\la}(c_{\la})$ and investigate the activity of critical points while we work with $\Bc=dd^c L(\la)$
and investigate the instability of cycles.\\
Using Branner-Hubbard holomorphic motion and quasi-conformal surgery, Uhre \cite{Uh} has built a motion similar to our in the region $U_1$.
The laminarity properties of the bifurcation current in the family of cubic polynomials have been studied
by Dujardin \cite{Du}, his approach differs from our.\\

The paper is organized as follows. The section \ref{premessa} describes the general concepts and results used in the paper, we discuss in particular the Milnor compactification 
of ${\cal M}_2$. The equidistribution formulas for $\Bc$ and $\Bc^k$ are
proved in section \ref{formulas} in the context of holomorphic families, these results remain valid for moduli spaces. In section \ref{HolMo},
we introduce the notion of guided holomorphic motion in ${\bf C}^2$ and establish, for these motions, an extension principle.
The section \ref{unif} is devoted to the uniformization of hyperbolic components in ${\cal M}_2$. 
 Finally, the laminarity properties of the bifurcation current in ${\cal M}_2$
are studied in sections \ref{Lamina} and \ref{U1}.

\section{The framework}\label{premessa}
\subsection{Bifurcation current and Lyapunov function}
 Every rational map of degree $d\ge 2$ on the Riemann sphere admits a maximal entropy measure
$\mu_f$. This measure satisfies $\pi^*\mu_f=dd^c G_F$
 where $\pi$ is the canonical projection from
${\bf C}^2\setminus\{0\}$ onto ${\bf P}^1$ and $G_F:=\lim d^{-n}\ln \Vert F^n\Vert$ is the Green function of any lift $F$ of $f$ to ${\bf C}^2$. 
The Lyapunov exponent of $f$ with respect to the measure $\mu_f$ is defined by $L(f)=\int_{{\bf P}^1} \ln \vert f'\vert\mu_{f}$
(see the book \cite{Sib} for a general exposition in any dimension).\\

When $f:M\times {\bf P}^1\to {\bf P}^1$ is a holomorphic family of degree $d$ rational maps, the Lyapunov function $L$ on the parameter space $M$
is defined by:
$$L(\la)=\int_{{\bf P}^1} \ln \vert f_{\la}'\vert\mu_{\la}$$
where $\mu_{\la}$ is the maximal entropy measure of $f_{\la}$. It turns out that the function $L$ is $p.s.h$ and H\"older continuous on $M$
(see \cite{BB} Corollary 3.4). The bifurcation current $\Bc$ of the family is a closed, positive $(1,1)$-current on $M$ which may be defined by
$$\Bc:=dd^c L(\la).$$
As it has been shown by DeMarco \cite{DeM2}, the support of $\Bc$ concides with the bifurcation locus of the family in the classical sense
of Ma\~n\'e-Sad-Sullivan (see also \cite{BB}, Theorem 5.2).
All these constructions also make sense on the moduli spaces ${\cal M}_d$ (see \cite{BB}, Section 6).

\subsection{The moduli spaces ${\cal M}_d$ and the hypersurfaces $Per_n(w)$}\label{M2}

The space $Rat_d$ of all rational maps of degree $d$ on
${\bf P}^1$ may be viewed as an open subset of ${\bf P}^{2d+1}$
on which the group of M\"obius transformations, which is isomorphic to 
$PSL(2,{\bf C})$, acts by conjugation.
The moduli space ${\cal M}_d$ is, by definition, the quotient resulting
from this action. We shall denote as follows the canonical
projection:

\begin{eqnarray*}
Rat_d &\longrightarrow& {\cal M}_d\\
f&\longmapsto&\bar f 
\end{eqnarray*}

Although the action of $PSL(2,{\bf C})$ is not free, it may be proven that 
${\cal M}_d$ is a normal quasi-projective variety \cite{Sil}.
For simplifying we shall sometimes commit the abuse of language which consists
in considering a parameter $\la \in {\cal M}_d$ as a rational map.
For instance, "$\la \in {\cal M}_d$ has a $n$-cycle" means that
every $f\in Rat_d$ such that $\bar f=\la$
posseses such a cycle.

We shall now recall why the subsets of ${\cal M}_d$ whose elements are the $\la$
which have a $n$-cycle of multiplier $w$ are (at least when $w\ne 1$)
hypersurfaces of ${\cal M}_d$.\\
Let $n\in \NN^*$. We first consider the hypersurfaces $Q_n^*$ of $Rat_d\times
{\bf P}^1$ defined by 
$$Q_n^*:=\{(f,z)\in Rat_d\times{\bf P}^1 \;/ \;f^n(z)=z\}$$
and then only retain the components whose generic points $(f,z)$ has the 
property that $z$ is exactly of period $n$:
$$Q_n:=\overline{Q_n^*\setminus\bigcup_{k<n, k\vert n}Q_k^*}.$$
We may now define a holomorphic function $w_n$ on $Q_n$ by setting
\begin{eqnarray*} 
w_n:Q_n&\longrightarrow& {\bf C}\\
  (f,z)&\longmapsto& (f^n)'(z)
\end{eqnarray*}
and consider its graph $\Gamma_n\subset Q_n\times{\bf C}\subset Rat_d\times {\bf P}^1
\times {\bf C}$. We are lead to the following key point:

\begin{defiprop}
Let $Per_n$ be the image of $\Gamma_n$ under the canonical projection
 $Rat_d\times {\bf P}^1
\times {\bf C}\to {\cal M}_d\times {\bf C}$. Then, $Per_n$ is an algebraic
subvariety of ${\cal M}_d\times {\bf C}$ and the projection
$Per_n\to{\cal M}_d$ is an analytic covering map. By construction, $(\bar f,w) 
\in Per_n$ if and only if $f$ there exists $z\in {\bf P}^1$ such that $(f,z)\in Q_n$
and $w_n(f,z)=w$.  
\end{defiprop}

We denote by $\nu_d(n)$ the number of distinct $n$-periodic points for a generic
$f\in Rat_d$. Then, $N(n):=\frac{\nu_d(n)}{n}$ is the cardinality of the 
generic fibers of the cover $Per_n\to {\cal M}_d$ and taking multiplicity into 
account, the fiber of any $\la \in {\cal M}_d$ may be written
$$w_{n,1}(\la),...,w_{n,N(n)}(\la).$$

We are now ready to give a precise description of the hypersurfaces $Per_n
(w)$

\begin{defiprop}\label{defifunct}
The hypersurface $Per_n\subset {\cal M}_d\times {\bf C}$ is defined by the 
algebraic equation $p_n(\la,w)=0$ where
$$p_n(\la,w):=\Pi_{j=1}^{N(n)}(w-w_{n,j}(\la)).$$
For any $w\in {\bf C}$, the hypersurface $Per_n(w)$ of ${\cal M}_d$ is defined
by:
$$Per_n(w):=\{\la\in {\cal M}_d /p_n(\la,w)=0\}.$$
If $w\ne 1$: $\la\in Per_n(w) \Leftrightarrow \la \;\textrm{posseses a cycle
of multiplier} \;w\;\textrm{and exact period}\; n $.
\end{defiprop}

The last assertion essentially follows from the fact that $(f,z)$ is a smooth 
point of $Q_n$ if $\bar f \in Per_n(w)$ when $w\ne 1$. 
This is a consequence of the implicit function theorem since, in that case,
$\frac{\partial}{\partial  z} (f^n(z)-z)=w_n(f,z)-1\ne 0$. Although much more 
delicate, a specific description of $Per_n(1)$ is possible
(see \cite{Mi} and \cite{Eps}). Let us stress that, in the present work, we 
will not have to use the 
hypersurfaces $Per_n(1)$.

\subsection{Milnor's compactification of ${\cal M}_2$}\label{compa}

We recall here some fundamental geometric properties
of the moduli space of quadratic rational maps. Our main references
for this subsection are the paper of Milnor \cite{Mi} and the fourth chapter of the
book of Silverman \cite{Silbook}. 

A generic $f\in Rat_2$ has $3$ fixed points with multipliers 
$\mu_1,\mu_2,\mu_3$. The symmetric functions 
$$\sigma_1:=\mu_1+\mu_2+\mu_3,\;\;\;\sigma_2:=\mu_1\mu_2+\mu_1\mu_3+\mu_2\mu_3,
\;\;\;\sigma_3:=\mu_1\mu_2\mu_3$$
are clearly well defined on ${\cal M}_2$ and it follows from the holomorphic index formula 
that $\sigma_3-\sigma_1 +2=0$. Milnor has actually shown that $(\sigma_1,\sigma_2)$ induces a good 
parametrization of ${\cal M}_2$ (\cite{Mi}).

\begin{theo}$({\bf Milnor})$\label{MilParam}
The map ${\cal M}_2\to{\bf C}^2$ defined by $\bar f \mapsto (\sigma_1,\sigma_2)$
is a biholomorphism. Using this identification, $p_n(\la,w)$ is a polynomial on 
${\bf C}^2\times {\bf C}$. Moreover, for every fixed $w\in {\bf C}$, the degree 
of $p_n(\cdot,w)$ is equal to $\frac{\nu_2(n)}{2}$ which is the number of 
hyperbolic components of period $n$ in the Mandelbrot set.
\end{theo}  

This gives a projective compactification
\begin{eqnarray*}
{\cal M}_2 \ni \bar f \longmapsto (\sigma_1:\sigma_2:1)\in {\bf P}^2
\end{eqnarray*}
whose corresponding line at infinity
will be denoted by ${\cal L}$
\begin{eqnarray*}
{\cal L}:=\{(\sigma_1:\sigma_2:0);\;(\sigma_1,\sigma_2)\in{\bf C}^2\setminus\{0\}\}.
\end{eqnarray*}

Any $Per_n(w)$ may be seen as a curve in ${\bf P}^2$.
It is important to stress that this compactification is actually natural and 
that the "behaviour near ${\cal L}$'' captures a lot of dynamically meaningful 
informations. We shall use the following facts, also due to Milnor 
(\cite{Mi}). 

\begin{prop}\label{Per/L}$({\bf Milnor})$
\begin{itemize}
\item[1)] For all $w\in{\bf C}$ the curve $Per_1(w)$ is a line  whose equation in ${\bf C}^2$ is
$(w^2+1)\la_1-w\la_2-(w^3+2)=0$ and whose
point at infinity is $(w:w^2+1:0)$. 
In particular, $Per_1(0)=\{\la_1=2\}$ is the line of 
quadratic polynomials, its point at infinity is $(0:1:0)$.
\item[2)] For $n>1$ and $w\in {\bf C}$ the points at infinity of the curves
$Per_n(w)$ are of the form $(u:u^2+1:0)$ with $u^q=1$ and $q\le n$.
\end{itemize}
\end{prop}

Working with the same compactification, Epstein \cite{Eps} has proved the boundedness of certain hyperbolic components of ${\cal M}_2$:

\begin{theo}\label{thEps}$({\bf Epstein})$ 
Let $H$ be a hyperbolic component of ${\cal M}_2$ whose elements are possessing two distinct attracting cycles. 
If neither attractor is a fixed point then $H$ is relatively compact in ${\cal M}_2$
\end{theo}

\section{Approximation formulas}\label{formulas}

In all the section we consider a fixed holomorphic family 
$f:M\times {\bf P}^1\to {\bf P}^1$ of degree $d$ rational maps 
whose parameter space $M$ is an $m$-dimensional complex manifold.\\

As for the moduli space ${\cal M}_d$ (see section \ref{premessa}), one associates to such a  
family a collection of polynomials
and hypersurfaces of $M$:

\begin{defiprop}\label{definova} Let $f:M\times {\bf P}^1\to {\bf P}^1$ be a family of degree $d$ rational maps on the Riemann sphere.
For any $n\in \NN^*$ there exists a monic polynomial $p_n\in {\cal O}(M)[X]$ such that for any
$w\in {\bf C}\setminus \{1\}$:
$$
 p_n(\la,w)=0 \Leftrightarrow
 f_{\la}\; 
\textrm{has a cycle of multiplier}\; w \;\textrm{and exact period}\; n.$$
The hypersurfaces $Per_n(w)$ are defined by $Per_n(w):=\{p_n(\cdot,w)=0\}$.   
\end{defiprop} 
We then define the currents $[Per_n(w)]$ by setting: 
$$[Per_n(w)]:=dd^c \ln \vert p_n(\cdot,w)\vert.$$
Up to multiplicity, $[Per_n(w)]$ is the current of integration on
the hypersurface $Per_n(w)$.\\
The main goal of this section is to show that the bifurcation current 
$\Bc$ of the family $(f_\la)_{\la\in M}$ is a limit of 
laminar currents whose dynamical interpretation is clear:

$$\Bc=\lim_n d^{-n}[Per_n(0)]=
\lim_n \frac{d^{-n}}{2\pi}\int_{0}^{2\pi}[Per_n(e^{i\theta})]d\theta.$$
Our proof essentially relies on an approximation formula for $L$
in terms of potentials of the form $\ln \vert p_n(\cdot,w)\vert$
(see theorem \ref{appL}) and on standard compacity properties for positive currents.\\

\subsection{Approximation formulas for the Lyapunov function}

\begin{defi}\label{d1}We use the notations introduced in the definition \ref{definova}.
For every integer $n$ we define a pair of functions $L_n$ and $L_n^0$ on the parameter space 
$M$ by setting:

$$L_n^0=d^{-n}\log \vert p_n(\la,0)\vert$$
$$L_n=\frac{ d^{-n}}{2\pi}\int_0^{2\pi}\log \vert p_n(\la,e^{i\theta})\vert d\theta.$$
\end{defi}

These functions are actually good approximations of the Lyapunov function $L$.
As we shall see, this basically follows from the fact that 
 the set of attracting cycles is finite and the following property of the Lyapunov function $L$:

\begin{theo}\label{thBDM}
Let $f$ be a rational map of degree $d$ and $RP(f,n)$ be the set of repelling $n$-periodic points of $f$. Then:
$$L(f)=\lim_n d^{-n}\sum_{RP(f,n)}\log\vert f'\vert=\lim_n \frac{d^{-n}}{n}\sum_{RP(f,n)}\log\vert (f^{n})'
\vert.$$
\end{theo} 

A similar property actually holds for holomorphic endomorphisms of ${\bf P}^k$ has it has been proved in \cite{BDM}.
The proof is more direct when $k=1$; basically one simply has to combine Lyubich's construction of repulsiv $n$-periodic points \cite{Lyu}
with Koebe theorem to control $(f^n)'$ on them. More precisely, while in the proof given in \cite{BDM} one uses a specific theorem
(the main result of that paper) to get the estimate $(a_2)$ that is to compare $\frac{1}{n}\ln \vert Jac f^n(z)\vert $ with $L(f)$,
one may simply use Cauchy estimates and the Koebe distorsion theorem.\\

Let us now state the approximation result on which will be based our study. 
 
\begin{theo}\label{appL}
We use the notations introduced in the above definition.
\begin{itemize}

\item[1)] The functions $L_n^0$ are $p.s.h$ on $M$.
The sequence $(L_n^0)_n$ is locally uniformly bounded from above and converges pointwise and in $L_{loc}^1$ to $L$.
\item[2)] The functions $L_n$ are positive, continuous and $p.s.h$ on $M$.
The sequence $(L_n)_n$ is locally uniformly bounded and converges pointwise (and therefore in $L_{loc}^1$) to $L$.  
\end{itemize}
\end{theo}

\proof
1) The function $L_n^0$ is clearly $p.s.h$ on $M$. 
 Let us denote by $w_{n,j}(\la)$ the roots of the polynomial $p_n(\la,\cdot)$:
$$p_n(\la,w)=\Pi_{j=1}^{N(n)}\left(w-w_{n,j}(\la)\right).$$
Using this decomposition one gets:

\begin{eqnarray}\label{i}
 L_n^0(\la)&=&d^{-n}\sum_{j=1}^{N(n)}{\log}\vert w_{n,j}(\la)\vert\le d^{-n}\sum_{j=1}^{N(n)}{\log}^+\vert w_{n,j}(\la)\vert=\nonumber\\
&=&\frac{d^{-n}}{n}\sum_{RP(f_{\la},n)}\log\vert (f_{\la}^{n})'\vert=d^{-n}\sum_{RP(f_{\la},n)}\log\vert f_{\la}'\vert.
\end{eqnarray}

As $f_{\la}^n$ has at most $d^n+1$ fixed points, this shows that the sequence $(L_n^0)_n$ is locally uniformly bounded from above. Moreover,
since $f_\la$ has only a finite number of attracting cycles and $w_{n,j}(\la)$ is the multiplier of such a cycle of period $n$
when $\vert w_{n,j}(\la)\vert<1$,
 the inequality in \ref{i} is an equality for all but a finite number of values of $n$. Then, for $n$ big enough we have:
\begin{eqnarray}\label{ii}
L_n^0(\la)= d^{-n}\sum_{j=1}^{N(n)}{\log}^+\vert w_{n,j}(\la)\vert=\frac{d^{-n}}{n}\sum_{RP(f_{\la},n)}\log\vert (f_{\la}^{n})'\vert.
\end{eqnarray}

 Thus, $(L_n^0)_n$ is a sequence of $p.s.h$ functions which
 is locally uniformly bounded from above and, according to theorem \ref{thBDM}, converges pointwise to $L$. 
By Slutsky's lemma, $L$ is therefore the only limit value of 
$(L_n^0)_n$ in $L_{loc}^1$. Then, by a classical compacity theorem for families of $p.s.h$ functions (see \cite{Sib} appendix),    
$(L_n^0)_n$ converges to $L$ in $L_{loc}^1$.

2)
Using the same decomposition than above for $p_n(\la,w)$ and the  
classical formula ${\log}^+\vert a\vert =\frac{1}{2\pi}\int_0^{2\pi}\log \vert e^{i\theta}-a\vert d\theta$, one gets:
\begin{equation}\label{1}
L_n(\la)=d^{-n}\sum_{j=1}^{N(n)}{\log}^+\vert w_{n,j}(\la)\vert.
\end{equation}
As $p_n(\la,w)$ is holomorphic in $\la$, this
shows that the function $L_n$ is positive, continuous and $p.s.h$ on $M$.
The sequence $(L_n)_n$ is, as before, locally uniformly bounded from above.
It now suffices to observe that $L_n^0(\la)=L_n(\la)$ for $n$ big enough, as it follows from \ref{1} and \ref{ii}.
\qed

\subsection{The bifurcation current as a limit of uniformly laminar currents}

Since the bifurcation current $\Bc$ is given by $\Bc=dd^cL$, our approximation formulas for $L$ lead to the following result:

\begin{theo}\label{appT} We use the notations introduced in the definition \ref{definova}. 
The bifurcation current $\Bc$ of a holomorphic family of rational maps coincides with the limit of the following sequences of laminar currents:
\begin{itemize}
\item[1)] $\Bc=\lim_n d^{-n}[Per_n(0)]$
\item[2)] $\Bc= \lim_n \frac{d^{-n}}{2\pi}\int_{0}^{2\pi}[Per_n(e^{i\theta})]d\theta.$
\end{itemize}
\end{theo}

The proof follows immediately from theorem \ref{appL}
by taking $dd^c$. For the second assertion however, one first has to observe that
the following identity occurs:

\begin{lem}\label{LemL_n}
$dd^c L_n=\frac{d^{-n}}{2\pi}\int_{0}^{2\pi}[Per_n(e^{i\theta})]d\theta.$
\end{lem}

\proof
Let  $\phi$ be a $(m-1,m-1)$ test form on $M$ where $m$ is the complex dimension of $M$.
One has to check that:
\begin{equation}\label{2}
\langle dd^c L_n,\phi\rangle=\frac{d^{-n}}{2\pi}\int_{0}^{2\pi} \langle[Per_n(e^{i\theta})],\phi\rangle d\theta.
\end{equation}
Let $dV$ be a volume form on $M$ and $\varphi$ be a smooth function such that $dd^c\phi=\varphi dV.$ 
Let us temporarily assume that the function
$\varphi(\la)\log\vert p_n(\la,e^{i\theta})\vert$ is integrable on $M\times[0,2\pi]$. Then \eqref{2} easily follows from Fubini's theorem: 

\begin{eqnarray*}
\langle dd^c L_n,\phi\rangle
&=&\langle L_n,dd^c\phi\rangle =\frac{d^{-n}}{2\pi}\int_{M}\big(\int_{0}^{2\pi}\log\vert p_n(\la,e^{i\theta})\vert d\theta\big) dd^c\phi\\
&=&\frac{d^{-n}}{2\pi}\int_{M}\big(\int_{0}^{2\pi}\varphi(\la)\log\vert p_n(\la,e^{i\theta})\vert d\theta\big) dV\\
&=&\frac{d^{-n}}{2\pi}\int_{0}^{2\pi}\big(\int_{M}\varphi(\la)\log\vert p_n(\la,e^{i\theta})\vert  dV\big)d\theta\\
&=&\frac{d^{-n}}{2\pi}\int_{0}^{2\pi}\langle dd^c \log\vert p_n(\la,e^{i\theta})\vert,\phi\rangle d\theta.\\ 
\end{eqnarray*}

It remains to see that $\log\vert p_n(\la,e^{i\theta})\vert$ is integrable on $\big(Supp\; \phi\big)\times[0,2\pi]$.
Let $c_n$ be an upper bound for $\log\vert p_n(\la,e^{i\theta})\vert$ on $\big(Supp\; \phi\big)\times[0,2\pi]$.
Then, the negative function   
$\log\vert p_n(\la,e^{i\theta})\vert -c_n$
is indeed integrable on $\big(Supp\; \phi\big)\times[0,2\pi]$.
This follows from the positivity of $L_n$ (see theorem \ref{appL}: 

\begin{eqnarray*}
\int_{Supp \; \phi}\big(\int_{0}^{2\pi}\big( \log\vert p_n(\la,e^{i\theta})\vert -c_n\big) d\theta\big)dV&=& 2\pi d^n\int_{Supp\; \phi}
L_n dV - 2\pi c_n \int_{Supp\; \phi}dV\\
&\ge& - 2\pi c_n \int_{Supp\; \phi}dV.\\
\end{eqnarray*}
\qed

Let us formulate some immediate consequences of the first assertion of theorem \ref{appT}.

\begin{cor}\label{AccCent}
A point $\la_0$ is in the bifurcation locus of the family
$f:M\times \pp\to\pp$ if and only if 
$\la_0=\lim_n \la_n$ 
where $f_{\la_n}$ has a super-attracting cycle of period $k_n$  
for some increasing sequence of integers $(k_n)_n$.
\end{cor}

\proof If $\la_0$ belongs to the bifurcation locus then the existence of the sequence $(\la_n)_n$ follows immediately
from theorem \ref{appT}. Conversely, the existence of such a sequence $(\la_n)_n$ implies that $
\la_0$ is not stable since, otherwise, $f_{\la_0}$ would have an infinite number of attracting cycles.\qed\\ 

Let us also recall that for the polynomial family $f_{\la}:=z^d+\la \;(\la\in{\bf C})$,
the bifurcation locus is the boundary of the Mandelbrot set ${ M}_d$ and the bifurcation
current is a measure which turns out to be the harmonic measure of ${M}_d$. In this setting the first assertion of 
theorem \ref{appT} yields the following corollary which was first proved by Levin \cite{Lev}:

\begin{cor}\label{mandel}
let ${M}_d$ be the set of complex numbers $\la$ for which the Julia set of $f_{\la}:=z^d+\la$ is connected. Then
$$\lim_n d^{-n}\sum_{f_{\la}^n(0)=0}\delta_{\la}=\mu$$
where $\mu$ is the harmonic measure of ${M}_d$.
\end{cor}

\subsection{Approximation of the bifurcation measure}

The higher degree 
 bifurcation currents $\big({\Bc}\big)^k$ ($k\le m$) of  holomorphic families 
have been studied in \cite{BB}. As one may expect, they can be obtained as limits of intersections of the laminar currents 
introduced in the former subsection.
For the sake of brevity we shall only consider here the case of the bifurcation measure
$\Bm=\left(dd^c L\right)^m$.\\

Our aim is to approximate $\Bm$ by measures of the following type:
\begin{equation}\label{3}
\int_{[0,2\pi]^m}[Per_{n_1}(e^{i\theta_1})]\wedge\cdot\cdot\cdot\wedge [Per_{n_m}(e^{i\theta_m})]\;d\theta_1\cdot\cdot\cdot d\theta_m.
\end{equation}
Let us observe that, for any fixed variety $Per_p(e^{i\theta_p})$,
the set of $\theta\in [0,2\pi]$ for which $Per_p(e^{i\theta_p})$
shares a non trivial component with $Per_m(e^{i\theta})$ for some $m\in \NN^*$
is at most countable. This follows from Fatou's theorem on the
finiteness of the set of non-repelling cycles. In particular, the wedge 
products $[Per_{n_1}(e^{i\theta_1})]\wedge\cdot\cdot\cdot\wedge [Per_{n_m}(e^{i\theta_m})]$
 make sense for almost every $(\theta_1,\cdot\cdot\cdot,
\theta_m)\in[0,2\pi]^m$ and the measures given in \ref{3} are well defined.
Our result may be stated as follows:

\begin{theo}\label{AppMu} We use the notations introduced in the definition \ref{definova}. Let $\Bm$ be the bifurcation measure 
of a holomorphic family $(f_{\la})_{\la\in M}$ of rational maps.
There exists increasing sequences of integers $k_2(n),...,k_m(n)$ such that:
$$\Bm=\lim_n \frac{d^{-(n+k_2(n)+\cdot\cdot\cdot+k_m(n))}}{(2\pi)^m}\int_{[0,2\pi]^m}
[Per_{n}(e^{i\theta_1})]\wedge\bigwedge_{j=2}^m[Per_{k_j(n)}(e^{i\theta_j})]\;d\theta_1\cdot\cdot\cdot d\theta_m.$$ 
\end{theo}

\proof
The problem is local and we may therefore replace $M$
by ${\bf C}^m$.
 We shall use the two following lemmas:

\begin{lem}\label{LemEgorov}
If $S_n\to (dd^c L)^p$ for some sequence $(S_n)_n$ of closed, positive $(p,p)$-currents on $M$ then  
 $dd^c L_{k(n)}\wedge S_n\to(dd^c L)^{p+1}$ for some increasing sequence of integers $k(n)$. 
\end{lem} 

\begin{lem}\label{LemwedgeL_n}
 $dd^c L_{n_1}\wedge\cdot\cdot\cdot\wedge dd^c L_{n_m}=
\frac{d^{-(n_1+\cdot\cdot\cdot+n_m)}}{(2\pi)^m}
\int_{[0,2\pi]^m}\bigwedge_{k=1}^{m}[Per_{n_k}(e^{i\theta_k})]d\theta_1\cdot\cdot\cdot
d\theta_m.$
\end{lem}

Using the first lemma inductively we obtain some increasing sequences of integers
$k_2(n),\cdot\cdot\cdot,k_{m-1}(n)$ such that 
$dd^c L_{k_p(n)}\wedge\cdot\cdot\cdot\wedge dd^c L_{k_2(n)}\wedge dd^c L_n \to \left(dd^c L\right)^{p+1}$
for any $2\le p\le m-1$. Then, lemma \ref{LemwedgeL_n} immediately yields to the desired result.\\

The proof of lemma \ref{LemwedgeL_n} is similar to that of lemma \ref{LemL_n} and we shall omit it. We now establish lemma \ref{LemEgorov}.
Let us denote by $s_n$ the trace measure of $S_n$, as $M$ has been identified with ${\bf C}^m$ this measure is given by 
$s_n:=S_n\wedge (dd^c \vert z\vert^2)^{m-p}$. Since $S_n$ is positive, $s_n$ is positive as well.
Let us consider the sequence $(u_k)_k$ defined by $u_k:=L_k-L$.
According to theorem \ref{appL}, $(u_k)_k$ converges pointwise to $0$ and is locally uniformly bounded (the function $L$ is continuous).
 
By Egorov theorem we may find subsets $E_n$ of $M$ such that $(u_k)_k$ is uniformly converging to $0$ on $E_n^c$
and $s_n(E_n)\le \frac{1}{n}$. We then pick an increasing sequence of integers $k(n)$ such that
$\Vert u_{k(n)}\Vert_{\infty,E_n^c}\le \frac{1}{n}$.\\

Writting
$dd^c L_k\wedge S_n - dd^c L\wedge (dd^c L)^p=
dd^c u_k \wedge S_n +  dd^c L \wedge(S_n-(dd^c L)^p)$,
we get the following estimate for any $m-(p+1)$ test form $\phi$ on $M$:

\begin{eqnarray*}
\vert\langle dd^c L_k \wedge S_n,\phi\rangle - \langle (dd^c L)^{p+1},\phi \rangle\vert
\le
\vert\langle dd^c u_k \wedge S_n,\phi\rangle\vert +\\
+ \vert\langle dd^c L \wedge S_n,\phi\rangle - \langle (dd^c L)^{p+1},\phi\rangle\vert.
\end{eqnarray*}
The second term of the above estimate tends to zero when $n$ tends to infinity and thus it remains to check that
 $\vert\langle dd^c u_{k(n)} \wedge S_n,\phi\rangle\vert$ tends to zero too. To this purpose we shall use the fact that the positive current $S_n$ may be considered as a
$(p,p)$ form whose coefficients are measures which are dominated by the trace measure $s_n$. Let $\chi$ be a cutoff function which is identically
equal to 1 on the support of $\phi$. Using the domination property of $s_n$ and the fact that $(u_{k(n)})_n$ is locally uniformly bounded we get:
 
\begin{eqnarray*}
\vert\langle dd^c u_{k(n)} \wedge S_n , \phi\rangle\vert&=&\vert \langle S_n , u_{k(n)}dd^c \phi \rangle\vert 
\le \textrm{cst}\; \Vert dd^c \phi\Vert_{\infty}\langle s_n , \vert u_{k(n)}\vert \chi\rangle \le\cdot\cdot\cdot\\
\cdot\cdot\cdot&\le& \textrm{cst}\; \Vert dd^c \phi\Vert_{\infty}\left(\langle s_n , \vert u_{k(n)}\vert \chi 1_{E_n}\rangle+\langle s_n , \vert u_{k(n)}\vert \chi 1_{E_n^c}\rangle\right)\\
&\le& \textrm{cst}\; \Vert dd^c \phi\Vert_{\infty} (\frac{1}{n} +\frac{1}{n}\langle s_n ,\chi\rangle).
\end{eqnarray*}

The conclusion follows since $s_n$ converges to the trace measure of $\left(dd^c L\right)^p$.
\qed\\

The theorem \ref{AppMu} has the following interesting consequence.

\begin{cor}
Let ${\cal N}_m$ be the set of parameters $\la$ such that $f_{\la}$ has $m$ distinct neutral cycles.
Then the support of the bifurcation measure $\Bm$ is contained in the closure of ${\cal N}_m$. 
\end{cor}

This corollary implies the existence of degree $d$ rational maps having $2d-2$ distinct neutral cycles,
we refer to \cite{BB} (Theorem 5.5 and Propositions 6.3, 6.8) for more details. This result was first proved by Shishikura
\cite{Shi} who also proved that any rational map of degree $d$ cannot have more than
 $2d-2$ distinct non-repulsiv cycles, his methods were based on quasi-conformal surgery.

\section{Holomorphic motions in ${\bf C}^2$}\label{HolMo}

This section is devoted to the study of a class of holomorphic motions in ${\bf C}^2$ which will appear naturally when we shall 
investigate the laminar structure of the bifurcation locus in the moduli space 
${\cal M}_2$. We essentially establish an extension principle which plays, for this class of motions,
the role played by the $\la$-lemma for holomorphic motions in ${\bf P}^1$ (see theorem \ref{thext}).

\subsection{Guided holomorphic motions}

\begin{defi}\label{defimot}
To any polynomial $p(\la,w):=\sum_{j=0}^{N} A_j(\la)w^j$ on ${\bf C}^2\times{\bf C}$ we associate the subset $\Omega_p$ of 
${\bf C}^2$ defined by:
$$\Omega_p:=\{\la\in {\bf C}^2\;\textrm{/}\;p(\la,\cdot) \;\textrm{has exactly one root in} \;\Delta\}.$$  
A \underline{$p$-guided holomorphic disc} is a holomorphic disc $\sigma:\Delta\to\overline{\Omega_p}$ such that 
$$p\big(\sigma(t),t\big)=0\;\;\textrm{for all}\;\; t\in \Delta.$$\\
 A \underline{$p$-guided holomorphic motion} ${\cal F}$ is a family of holomorphic discs which are $p$-guided and  mutually disjoint.
\end{defi}

 Let us emphasize that in the above definition the set $\Omega_p$ is not supposed to be open. Note also that 
$\la\notin\overline{\Omega_p}$ when $p(\la,\cdot)$ has at least two roots in $\Delta$.
 This yields to the following simple, but important, observation.
We recall that two holomorpic discs $\sigma_1,
\sigma_2:\Delta\to{\bf C}^2$ are said to be {\it glued} if there exists
 non-constant holomorphic functions $u_1,u_2:\Delta_r\to\Delta$ such that
$\sigma_1\circ u_1=\sigma_2\circ u_2$ on $\Delta_r$.

\begin{prop}\label{obs}
Two $p$-guided holomorphic discs $\sigma_1$ and $\sigma_2$ which are glued (or, in particular, have the same image)
must coincide: $\sigma_1(t)=\sigma_2(t)$ for all $t\in \Delta$. Every $p$-guided holomorphic disc $\sigma:\Delta\to 
\overline{\Omega_p}$ is one-to-one.
\end{prop}

\proof if $\sigma_2\circ u_2(t)=\sigma_1\circ u_1(t)=:\la_t\in\overline{\Omega_p}$ for $t\in \Delta_r$ then
$p(\la_t,u_1(t))=p(\la_t,u_2(t))=0$ and therefore $u_1\equiv u_2$ on $\Delta_r$. It follows that $\sigma_1$
and $\sigma_2$ coincide on some open subset of $\Delta$ and thus everywhere.\qed\\

It is useful to see the guided holomorphic motions as given by special kinds of mappings. To this purpose,
for any family of discs $\cal{F}$ and any $t_0\in \Delta$, we set $${\cal F}_{t_0}:=\{\sigma(t_0)\;;\; \sigma\in {\cal F}\}.$$
When ${\cal F}$ is a $p$-guided holomorphic motion and $\la\in{\cal F}_{t_0}$, 
there exists a unique holomorphic disc $\sigma_\la \in {\cal F}$
such that $\sigma_{\la}(t_0)=\la$.
 Thus, for any fixed $t_0\in \Delta$, the holomorphic motion ${\cal F}$ may be parametrized by a map 

\begin{eqnarray*}
\sigma:{\cal F}_{t_0}\times\Delta&\longrightarrow& \overline{\Omega_p}\\
(\la,t)&\longmapsto&\sigma_{\la}(t)
\end{eqnarray*}

satisfying the following properties:

\begin{itemize}
\item[1)]$\sigma_{\la}(\cdot)$ is holomorphic on $\Delta$ for every $\la\in{\cal F}_{t_0}$
\item[2)]$\sigma_{\la}(t_0)=\la$ for every $\la\in{\cal F}_{t_0}$
\item[3)]$\sigma_{\la}(\Delta)\cap\sigma_{\la'}(\Delta)=\emptyset$ if $\la\ne\la'$
\item[4)]$p\big(\sigma_{\la}(t),t\big)=0$ for every $\la\in{\cal F}_{t_0}$ and every $t\in \Delta$.
\end{itemize}

Let us underline that our notations might be confusing if compared with those usually 
adopted for holomorphic motions on the Riemann sphere: for us $\la$ is not the ``holomorphic time''
but a point under motion !\\ 
We will often identify a $p$-guided holomorphic motion with one of its parametrizations 
$\sigma:{\cal F}_{t_0}\times\Delta\longrightarrow \overline{\Omega_p}$. In particular, this point of view
allows to define in a simple way the notion of {\it continuous $p$-guided holomorphic motion}. 
The property of continuity will play an important role in our study.

\begin{defi}
A $p$-guided holomorphic motion $\cal F$ is said to be continuous if for any $t_0\in \Delta$
the parametrization
$\sigma:{\cal F}_{t_0}\times\Delta
\longrightarrow \overline{\Omega_p}$ is continuous.
\end{defi}

We may now state the main result of the section.

\begin{theo}\label{thext}
Let $p(\la,w)$ be a polynomial on ${\bf C}^2\times{\bf C}$ such that the degree of $p(\cdot,w)$ does not depend on $w\in \Delta$.
Let ${\cal G}$ be a $p$-guided holomorphic motion in ${\bf C}^2$ such that any component of the algebraic curve $\{p(\cdot,t)=0\}$
contains at least three points of ${\cal G}_t$ for every $t\in \Delta$.
Then, for any ${\cal F}\subset{\cal G}$ such that ${\cal F}_{t_0}$ is relatively compact in ${\bf C}^2$ for some $t_0\in \Delta$,
there exists a continuous $p$-guided holomorphic motion $\widehat{{\cal F}}$ in ${\bf C}^2$
such that ${\cal F}\subset \widehat{{\cal F}}$
and $\widehat{{\cal F}}_{t_0}=\overline{{\cal F}_{t_0}}.$
\end{theo}
 
Let us emphasize that some hyperbolic components in the moduli space ${\cal M}_2$ 
will appear to be moved by holomorphic motions satisfying the assumptions of the above theorem (see subsection \ref{HM}).
The remaining of the section will be devoted to the proof of theorem \ref{thext}.

\subsection{A compactness property}

To start with, we summarize a few basic intersection properties of (guided) holomorphic discs.\\

\begin{prop}\label{propdiscs}
\begin{itemize}
\item[1)] Let $\sigma_n:
\Delta\to{\bf C}^2$ and $s_n:
\Delta\to{\bf C}^2$  
be two sequences of holomorphic discs which converge to 
$\sigma,
s:\Delta\to{\bf C}^2$. If $\sigma$ and $s$ are not glued and $\sigma(0)=s(0)$ then
$\sigma_n(\Delta)\cap s_n(\Delta)\ne \emptyset$ for $n$ big enough. 
\item[2)]Let $\sigma_n
:\Delta\to \overline{\Omega_p}$ and $s_n
:\Delta\to \overline{\Omega_p}$ 
 be two sequences of $p$-guided holomorphic discs such that   
$\sigma_n(\Delta)\cap s_n(\Delta)= \emptyset$ for all $n\in \NN$. Then the limits of $(\sigma_n)_n$ and $(s_n)_n$ either 
coincide or are disjoint.
\end{itemize}
\end{prop}

\proof The first assertion is classical, we include here a proof for the sake of completness.
We may assume that $\sigma(0)=s(0)=(0,0)$. Let $E$ be the analytic set defined by 
$E=\{(x,y)\in\Delta\times\Delta\;\textrm{/}\; \sigma(x)=s(y)\}$.
The point $(0,0)$ belongs to $E$ and, since the discs are not glued, is isolated in $E$. In other words, 
there exists $0<\rho<\frac{1}{2}$ such that
 $x=y=0 $ as soon as $\sigma(x)=s(y)$ and $\max(\vert x\vert,\vert y\vert)< 2\rho$.  
Let us consider the holomorphic maps $\varphi,\varphi_n: \Delta_{\rho}\times\Delta_{\rho}\to{\bf C}^2$ respectively defined by
$\varphi(x,y)=\sigma(x+y)-s(x-y)$ and $\varphi_n(x,y)=\sigma_n(x+y)-s_n(x-y)$. As $\varphi^{-1}\{(0,0)\}=\{(0,0)\}$ we must have 
$\varphi_n^{-1}\{(0,0)\}\ne \emptyset$ and in particular $\sigma_n(\Delta)\cap s_n(\Delta)\ne \emptyset$ for $n$ big enough 
(see for instance
\cite{BL} lemme 3.4).\\
To deduce the second assertion from the first, we observe that the limits of $p$-guided holomorphic discs are $p$-guided
and recall that 
two $p$-guided holomorphic discs which are glued must coincide (Prop. \ref{obs}).\qed\\

When proving theorem \ref{thext}, we shall need a substitute to the Montel-Picard theorem.
The following lemma will play this role, its proof is based on the Zalcman renormalization principle
(see \cite{Z} or \cite{Be}).

\begin{lem}\label{lemnorma}
Let $p(\la,w)$ be a polynomial on ${\bf C}^2\times{\bf C}$ and 
${\cal G}$ be a $p$-guided holomorphic motion.
Assume that both $p$ and 
$\cal G$  satisfy the assumptions of theorem \ref{thext}. Let
${\cal F}\subset {\cal G}$. If there exists a compact subset $K$ of ${\bf C}^2$ and $r\in [0,1[$ such that $\sigma(\Delta_r)\cap K\ne \emptyset$ for every 
$\sigma\in {\cal F}$, then ${\cal F}$ is a normal family in ${\cal O}\big(\Delta,{\bf C}^2\big)$.
\end{lem}

\proof
Since ${\cal F}\subset{\cal O}\big(\Delta,{\bf C}^2\big)$ and $\sigma(\Delta_r)\cap K$ is never empy when $\sigma \in {\cal F}$, it suffices to show that 
${\cal F}$ is normal in ${\cal O}\big(\Delta,{\bf P}^2\big)$. We proceed by contradiction and assume that the family ${\cal F}$ is not normal in 
${\cal O}\big(\Delta,{\bf P}^2\big)$. Then, according to Zalcman's renormalization principle, 
there exists $\{s_n\;;\;n\in \NN\}\subset {\cal F}$,
$t_n\in \Delta$ and $\rho_n>0$ such that  $\
lim_n\rho_n=0$, $\lim_nt_n=t_0$ for some $t_0\in \Delta$ and $s_n(t_n+\rho_n t)
\to h$. Moreover, the convergence of $s_n(t_n+\rho_n t)$ is locally uniform on 
${\bf C}$ and the limit $h$ is a non constant entire curve.\\
Let us show that $h({\bf C})\subset {\bf C}^2\cap {\cal P}_{t_0}$ where ${\cal P}_t$ denotes the closure in ${\bf P}^2$ of the curve 
$\{p(\cdot,t)=0\}$. Since
the motion 
is
guided we have $s_n(t_n+\rho_n t)\in {\cal P}_{t_n+\rho_n t}$. As the degree of $p(\la,t)$ in $\la$ does not depend on 
$t\in \Delta$, the family $\big({\cal P}_t\big)_{t\in \Delta}$ is a continuous family of curves in ${\bf P}^2$. 
One sees therefore, since $t_n+\rho_nt\to t_0$, that $h(t)\in {\cal P}_{t_0}$ for every $t\in{\bf C}$.
In particular, $h({\bf C})$ is not contained in the line at infinity and, since $s_n(\Delta )\subset{\bf C}^2$, the curve $h({\bf C})$ is actually contained in
${\bf C}^2$. We have shown that $h({\bf C})\subset {\bf C}^2\cap {\cal P}^o_{t_0}$ for some component ${\cal P}^o_{t_0}$ of ${\cal P}  _{t_0}$.\\
It remains to show that such an entire curve does not exist. Picard's first theorem implies that any non-constant entire curve in 
${\cal P}^o_{t_0}$ cannot avoid more than two points. Thus, since $h({\bf C})\subset {\bf C}^2\cap {\cal P}^o_{t_0}$ and $Card\big({\cal G}_{t_0}\cap  
 {\cal P}^o_{t_0}\cap{\bf C}^2\big)\ge 3$, the conclusion will follow if we show that $h({\bf C})\cap {\cal G}_{t_0}=\emptyset$. 
We shall actually see that $h({\bf C})$ does avoid every disc of ${\cal G}$.
Let us first observe that for $\sigma\in{\cal G}$ and $R>0$, the discs 
$\sigma:\Delta\to\overline{\Omega_p}$ and $h:\Delta_R\to{\cal P}^o_{t_0}\cap{\bf C}^2$ are not glued.
If they would be glued then 
$\sigma(\Delta)$ would be contained in ${\cal P}_{t_0}\cap{\bf C}^2$ and $p(\sigma(t),\cdot)$ would have both $t$ and $t_0$ as roots for any $t\in\Delta$
which is impossible for $\sigma(t)\in\overline{\Omega_p}$.
Then, as the discs of ${\cal G}$ are mutually disjoint, we have 
$\sigma(\Delta)\cap s_n(t_n+\rho_n R\Delta)=\emptyset$, and the first assertion of proposition \ref 
{propdiscs} shows indeed that $\sigma(\Delta)\cap h(\Delta_R)=\emptyset$.\qed

\subsection{Extension of guided holomorphic motions}

We prove the theorem \ref{thext}.\\

\proof \underline{Ste}p\underline{ 1}: {\it For any $t_0'\in \Delta$ and any $\la_0\in\overline{{\cal F}_{t_0'}}$ there exists a unique $p$-guided holomorphic disc 
$\sigma:\Delta\to\overline{\Omega_p}$
such that $\sigma_\la$ tends to $\sigma$ when $\la$ tends to $\la_0$ in ${\cal F}_{t_0'}$.}\\

According to lemma \ref{lemnorma}, $\cal F$ is normal in ${\cal O}\big(\Delta,{\bf C}^2\big)$. 
It therefore suffices to show that $\big(\sigma_{\la_n}\big)_n$ 
has only one limit
in ${\cal O}\big(\Delta,{\bf C}^2\big)$ when $\la_n$ tends to $\la_0$ in ${\cal F}_{t_0'}$. As every limit of 
$\big(\sigma_{\la_n}\big)_n$ goes through
$\la_0=\lim_n\la_n=\lim_n\sigma_{\la_n}(t_0')$, this follows immediately from the second assertion of proposition \ref{propdiscs}.\\

\underline{Ste}p\underline{ 2}: {\it Construction of $\widehat{\cal F}$.}\\

According to the first step (applied with $t_0'=t_0$),
 it makes sense to define a $p$-guided holomorphic disc $\sigma_{\la_0}$ for every ${\la_0}\in \overline{{\cal F}_{t_0}}$ by 
setting:
$$\sigma_{\la_0}:=\lim_{{\cal F}_{t_0}\ni\la\to\la_0}\sigma_\la.$$
Then, using the second assertion of proposition \ref{propdiscs}, one sees that the family 
$$\widehat{\cal F}:=\{\sigma_\la,\;\la\in\overline{{\cal F}_{t_0}}\}$$ consists in mutually disjoint discs and 
is therefore a $p$-guided holomorphic motion in ${\bf C}^2$. By construction, $\cal F\subset \widehat{\cal F}$ and 
$\widehat{{\cal F}}_{t_0}=\overline{{\cal F}_{t_0}}.$\\

\underline{Ste}p\underline{ 3}: {\it Continuity of $\widehat{\cal F}$.}\\

Let us set $\widehat{\cal G}:={\cal G}\cup\widehat{\cal F}$. Using again the second assertion of proposition \ref{propdiscs}, 
one sees that the discs of       
$\widehat{\cal G}$ are mutually disjoint. Thus $\widehat{\cal G}$ is a $p$-guided holomorphic motion in ${\bf C}^2$. According to the second step
 $\widehat{{\cal F}}_{t_0}=\overline{{\cal F}_{t_0}}$ is compact in ${\bf C}^2$ and then, as ${\cal G}\subset\widehat{\cal G}$,
 we may apply the fisrt step to the pair 
$(\widehat{\cal F},\widehat{\cal G})$. Thus, there exists a unique disc $\sigma$ such that,
for any $\la_0 \in \widehat{{\cal F}}_{t_0'}$ and any sequence $(\la_n)_n$ which tends to $\la_0$ in $\widehat{{\cal F}}_{t_0'}$, 
we have $\lim_n\sigma_{\la_n}=\sigma$.
Considering the stationnary sequence $\la_n=\la_0$ one gets that $\sigma=\sigma_{\la_0}$. Thus, $\lim_n\sigma_{\la_n}(t_n)=
\sigma(t_0')=\sigma_{\la_0}(t_0')$ when 
$(\la_n,t_n)$ tends to $(\la_0,t_0')$ in $\widehat{{\cal F}}_{t_0'}\times \Delta$. \qed

\section{Uniformization of some hyperbolic components in ${\cal M}_2$}\label{unif}
 
We identify ${\cal M}_2$ with ${\bf C}^2$ (see theorem \ref{MilParam}) and, for any $n\in \NN^*$, we consider the following open region:
\begin{eqnarray*}U_n:&=&\{\la\in{\bf C}^2/ f_{\la} \;\textrm{has an attracting cycle of period }n\}\\
&=&\{\la\in{\bf C}^2/ p_n(\la,\cdot) \;\textrm{has a root in }\Delta\}
\end{eqnarray*}
we  recall that the polynomials $p_n$ and the hypersurfaces $Per_n$ have been introduced in the section \ref{premessa}.\\
When $n\ne m$, the open set $H_{n,m}:=U_n\cap U_m$ consists of hyperbolic parameters $\la$ possessing two distinct attracting cycles.
The periods of these cycles are respectively equal to $n$ and $m$ while their multipliers are roots of $p_n(\la,\cdot)$ and $p_m(\la,\cdot)$. Since 
a quadratic rational map have at most two attracting cycles, the polynomials $p_n(\la,\cdot)$ and $p_m(\la,\cdot)$ have both exactly one root in 
$\Delta$ when $\la\in 
H_{n,m}$ and thus, denoting by $w_n(\la)$ and $w_m(\la)$ these roots, one sees that 
 the map $\la\mapsto\left(w_n(\la),w_m(\la) \right)$ is well defined and holomorphic on $H_{n,m}$. Let us formalize this.

\begin{defi}
For $n\ne m $ there exists a holomorphic map $\phi_{n,m}:H_{n,m}\to\Delta^2$ defined on the open set 
$H_{m,n}:=U_n\cap U_m$ by $\phi_{n,m}=\left(w_n(\la),w_m(\la) \right)$ where
$p_n(\la,w_n(\la))=p_m(\la,w_m(\la))=0$.
\end{defi}

A classical result due to Douady and Hubbard (\cite{DH},\cite{CG}) asserts that the bounded hyperbolic components of the polynomial quadratic family
are uniformized by the multiplier of the attracting cycle. Their proof is based on a quasi-conformal surgery argument. 
We will adapt it for showing that $\phi_{n,m}$ induces a biholomorphism between 
any connected component of $H_{n,m}$ and the bidisc $\Delta^2$.

\begin{theo}\label{thUni}
Let $n$ and $m$ be distinct positive integers. The map $\phi_{n,m}$ induces a biholomorphism $H_{n,m}^j\stackrel {\phi_{n,m}}\longrightarrow\Delta^2$
on any connected component $H^j_{n,m}$ of $H_{n,m}$.
\end{theo}

\proof
Since $\Delta^2$ is simply connected it suffices to show that $\phi_{n,m}$ is proper and locally invertible.\\

\underline{Ste}p\underline{ 1}: {\it Properness}.\\

The properness of $\phi_{n,m}$ will be deduced from Milnor's result (see Proposition \ref{Per/L}) when $n=1$ and from Epstein's boundedness result
 (see theorem \ref{thEps})
in the general case.\\
Assume that  
$H_{n,m}^j\stackrel {\phi_{n,m}}\longrightarrow\Delta^2$ is not proper. Then there exists a sequence $(\la_k)_k$ in 
$H_{n,m}^j$ such that
$(w_n(\la_k),w_m(\la_k))=\phi_{n,m}(\la_k)$ tends to $(\alpha,\beta)\in\Delta^2$ and which tends to the boundary of $ H_{n,m}^j$ or to infinity.
Let us first consider the case where both $n$ and $m$ differ from 1. Then, according to theorem \ref{thEps}, we may assume that
$\la_k\to \la_0\in bH_{n,m}^j$. Taking the limit in $p_n(\la_k,w_n(\la_k))=p_m(\la_k,w_m(\la_k))=0$ one gets
$p_n(\la_0,\alpha)=p_m(\la_0,\beta)=0$, which means that $\la_0\in U_n\cap U_m$ and is absurd.\\
Let us now discuss the case $n=1<m$. If the sequence $(\la_k)_k$ would not go to infinity we could simply argue as before. We 
will thus assume that $(\la_k)_k$ is converging to some $\la_0\in{\cal L}$ in ${\bf P}^2$.
Since $\la_k\in {Per}_1\left(w_1(\la_k)\right)$, the first assertion of Proposition \ref{Per/L} shows that $\la_0=[\alpha:\alpha^2+1:0]$
(we recall that the curves $Per_1(w)$ are lines in ${\bf C}^2$). On the other hand, as the degree in $\la$ 
of $p_m(\la,w)$ is constant when $w\in\Delta$, the family of curves $\overline {Per}_m\left(w\right)_{w\in\Delta}$ is continuous and therefore
$\la_0\in{\cal L}\cap\overline {Per}_m\left(\beta\right)$. By the second assertion of Proposition \ref{Per/L}, 
this implies that $\la_0=[e^{i\theta}:e^{2i\theta}:0]$ for some $\theta \in [0,2\pi[$ which, again,
is absurd.\\
 
\underline{Ste}p\underline{ 2}: {\it Local invertibility}.\\

The key point is the following lemma. Its proof will be sketched in the third step.

\begin{lem}\label{key}
Let $\la_0\in H_{n,m}^j$ and $(a,b):=\phi_{n,m}(\la_0)$. Then there exists a continuous map $t\mapsto ({\bar f}_t)$ from $\Delta_r$
($\vert b \vert <r$) to $ H_{n,m}^j\cap Per_n(a)$ such that ${\bar f}_b=\la_0$ and $w_m({\bar f}_t)=t$.
\end{lem}

It suffices to show that the curves $\{w_n=cst\}$ are smooth in $H_{n,m}^j$ and that $dw_m$ does not vanish on them. Indeed, by symmetry, the same will occur
for the curves $\{w_m=cst\}$ and $dw_n$, which easily implies that $dw_m$ and $dw_n$ are independant on $H_{n,m}^j$.\\
Let $a\in \Delta$. We first want to see that the curve $\{w_n=a\}=\{p_n(\cdot,a)=0\}$ is smooth in $H_{n,m}^j$. Assume to the contrary 
that this would not be the case at some point $\la_0\in H_{n,m}^j$. Then the singularity of $\{p_n(\cdot,a)=0\}$ at $\la_0$ is not a point of 
self-intersection
 since otherwise, for $a'$ close to $a$, 
$\{p_n(\cdot,a')=0\}$ would intersect $\{p_n(\cdot,a)=0\}$ at some point having three distinct attracting cycles (two of period $n$ and one of period $m$). 
Thus, the curve  $\{p_n(\cdot,a)=0\}$ has a cusp at $\la_0$ and therefore, the restriction of the holomorphic function $w_m$ to the curve is not one-to-one 
on 
any neighbourhood of $\la_0$. Composing with a local uniformization we get a holomorphic function $\tilde {w}_m$ on the unit disc which is not one-to-one
and, by the lemma \ref{key}, a continuous map $t\mapsto {\tilde f}_t$ from a neighbourhood of $b\in \Delta_r$ to the unit disc such that
${\tilde w}_m\left({\tilde f}_t\right)=t$: this is impossible.\\
Now, since $\{p_n(\cdot,a)=0\}$ is smooth at $\la_0$, we may repeat the last argument to see that the restriction of $w_m$ to the curve is one-to-one
near $\la_0$. This implies that $dw_m$ does not vanish on the tangent plane to the curve at $\la_0$.\\

\underline{Ste}p\underline{ 3}: {\it Proof of the lemma}.\\

We follow the proof of Douady-Hubbard as it is presented in \cite{CG} pages 134-135 and simply mention the points which need to be adapted.\\
Let us first justify that, for any $\la\in H_{n,m}^j$, the immediate basins of both attractors are simply connected. As the map  
$\phi_{n,m}:H_{n,m}^j\longrightarrow\Delta^2$ is proper holomorphic, it is onto and in particular $\phi_{n,m}^j(\la_c)=(0,0)$ for some
$\la_c\in H_{n,m}^j$. The attracting basins of $\la_c$ are two distinct super-attracting basins, each of them is simply-connected
since it contains only one of the two critical points of $\la_c$. This property is kept on $H_{n,m}^j$ because $H_{n,m}^j$ is a stable component.\\
We now construct the map $t\mapsto f_t\in Rat_2$. Let $f\in Rat_2$ be a representant of $\la_0$.
 After a change of coordinates we may assume that $\infty$ belongs to the $n$ cycle of $f$.
We perform exactly the same quasi-conformal surgery than Douady-Hubbard. Here however, instead of producing a quadratic polynomial we
want to have $w_n(f_t)=a$. In other words, all $f_t$ must have a $n$-cycle whose multiplier equals $a$.
Actually, the normalization at $\infty$ made in Douady-Hubbard's proof guarantees that. Indeed, $f_t=\psi_t\circ g_t \circ \psi_t^{-1}$ where $g_t$ coincides 
with $f$
near $\infty$ and $\psi_t$ is (holomorphic) tangent to identity at $\infty$.
Finally, the conditions $w_m(f_t)=t$ and $f_b=f$ as well as the continuity of $
f_t$ are the results of the surgery itself.\qed\\

We would like to mention that the local invertibility of the map $\phi_{n,m}$
may also be directly obtained from a recent general transversality result
whose proof is based on the quadratic differentials techniques \cite{BEH}: 

\begin{theo}\label{thBEH}$({\bf Epstein})$\\
For every $\eta, \eta' \in \Delta$ and every distinct $n,n' \in \NN^*$ the curves $Per_n(\eta)$ and $Per_{n'}(\eta')$ are 
smooth and intersect transversally at any common point.
\end{theo}

The following proposition shows that the curves $Per_n(w)$ have no multiplicity, it also
relates the degree of the polynomials $p_n(\cdot,w)$ with the number of hyperbolic components of period $n$ in the Mandelbrot set.

\begin{prop}\label{multi}
Let $N_2(n):=Card\left(Per_n(0)\cap Per_1(0)\right)$ be the number of hyperbolic components of period $n$ in the Mandelbrot set. Then $N_2(n)=\frac{\nu_2(n)}{2}$
where $\nu_2(n)$ is defined inductively by $\nu_2(1)=2$ and $2^n=\sum_{k\vert n}\nu_2(k)$.
Moreover, for any $w\in \Delta$ and any $\eta\in\Delta$ we have $Deg\;p_n(\cdot,w)=N_2(n)=Card\left(Per_n(w)\cap Per_1(\eta)\right).$
\end{prop}

\proof
Let us first compute $N_2(n)$. We set $P_n(c):=p_c^n(0)$ where $p_c(z)=z^2+c$. Then, $P_n$ is a polynomial of degree $2^{n-1}$ whose roots are simple
(see \cite{DH}
page
108). Since $P_n(c)$ vanishes if and only if the critical point $0$ is a periodic point for $p_c$ whose period is dividing $n$, we have
$N_2(n)=2^{n-1}-\sum_{k\vert n,k<n} N_2(k)$.\\
Let $w,\eta\in\Delta$. The existence of the biholomorphisms $\phi_{1,n}^j:H_{1,n}^j\longrightarrow\Delta^2$
(see theorem \ref{thUni}) implies that $Card\left(Per_n(w)\cap Per_1(\eta)\right)=Card\left(Per_n(0)\cap Per_1(0)\right)=N_2(n)$.
According to the proposition \ref{Per/L}, the line $Per_1(0)$ does not meet $Per_n(w)$ at infinity and thus
$Deg\;p_n(\cdot,w)\ge Card\left(Per_n(w)\cap Per_1(\eta)\right)=N_2(n)$. The function $Deg\;p_n(\cdot,w)$ being $l.s.c$ on $\Delta$, it remains to see that 
the above inequality is actually an equality on a dense subset of $\Delta$. This follows from the fact that $p_n$ is a defining function of $Per_n$ (see
the proposition \ref{defifunct}).\qed

\section{Laminar structures in the bifurcation locus of ${\cal M}_2$}\label{Lamina}

The main goal of this section is to show that the bifurcation current $\Bc$ is uniformly laminar in the following open regions of ${\cal M}_2$:
\begin{eqnarray*}
U_n:&=&\{\la\in{\bf C}^2/ f_{\la} \;\textrm{has an attracting cycle of period }n\}\\
&=&\{\la\in{\bf C}^2/ p_n(\la,\cdot) \;\textrm{has a root in }\Delta\}.
\end{eqnarray*}

We shall in particular prove that the bifurcation loci in the regions $U_n$ are
 obtained by moving holomorphically the bifurcation loci of their ``central curves''
$Per_n(0)$. As in the former section we set $H_{n,m}:=U_n\cap U_m$ and denote by $H_{n,m}^j$ the connected component of $H_{n,m}$.

\subsection{Holomorphic motion of centers}\label{HM}

We work here in a fix open region $U_n$. What we have called the {\it central curve of} $U_n$ is $Per_n(0)$.
The centers of the hyperbolic components $H_{n,m}^j$, that is the points
of $Per_n(0)\cap Per_m(0)$ for $m\ne n$,
 will be called the {\it centers of} $Per_n(0)$. To keep an example in mind, let us recall
that $Per_1(0)$ is the line of quadratic polynomials and its centers are exactly those of the hyperbolic components (except the main cardioid)
of the polynomial quadratic family.
The uniformization discussed in section \ref{unif} (theorem \ref{thUni}) naturally induces a holomorphic motion of the centers. Let us state this more formaly:

\begin{defiprop}\label{Mouvcenter}
A point $\la\in Per_n(0)$ is called a $n$-center if the corresponding rational maps have a super-attracting cycle of period $m\ne n$.
There exists a $p_n$-guided holomorphic motion ${\cal C}_n$ whose discs are contained in $\Omega_{p_n}$ and such that $\left({\cal C}_n\right)_0$
is exactly the set of $n$-centers. Moreover, ${\cal C}_n=\cup_{m\ne n}Per_m(0)\cap U_n$.  
\end{defiprop}

We recall that
this means that there exists a map 

\begin{eqnarray*}
\sigma:\left({\cal C}_n\right)_0\times\Delta & \longrightarrow &\Omega_{p_n}\subset U_n\\
(\la,t)&\longmapsto&\sigma_\la(t)\\
\end{eqnarray*}

which is holomorphic in $t$ and satisfies:

\begin{itemize}

\item[1)]$\sigma_{\la}(0)=\la$ \;$\forall \la\in\left({\cal C}_n\right)_0$
\item[2)]$\sigma_{\la}(\Delta)\cap\sigma_{\la'}(\Delta)=\emptyset$ if $\la\ne\la'$
\item[3)]$p_n\big(\sigma_{\la}(t),t\big)=0$ \; $\forall\la\in\left({\cal C}_n\right)_0,\;\forall t\in \Delta$. 
\end{itemize}

 \proof
We simply have to exhibit, for any $n$-center $\la_0$, a $p_n$-guided holomorphic disc $\sigma:\Delta\to\Omega_{p_n}$ such that
$\sigma(0)=\la_0$ (see definition \ref{defimot}).
By assumption, $\la_0$ belongs to some component $H_{n,m}^j$. We set $\sigma(t)=(\Phi_{n,m}^j)^{-1}(t,0)$ where $\phi_{n,m}^j$ is the biholomorphism
between $H_{n,m}^j$ and $\Delta^2$ given by theorem \ref{thUni}. By construction $\sigma(\Delta)\subset H_{n,m}^j$
and $p_n(\sigma(t),t)=0$. As a rational map of degree two has at most two attracting cycles, we have $\sigma(\Delta)\subset H_{n,m}^j\subset\Omega_{p_n}$.
The discs of ${\cal C}_n$ are mutually disjoint since they belong to distinct connected components.\qed\\

The following lemma is essentially based on the first assertion of theorem \ref{appT}, it will be used to transmit the motion from centers to the bifurcation locus.
Let us recall that $\Bc$ denotes the bifurcation current in ${\cal M}_2$.

\begin{lem}\label{lemcent}
\begin{itemize}

\item[1)] If $\la_0 \in U_n\cap \textrm{Supp}(\Bc)$ then there exists $a_m\in ({\cal C}_n)_0$, $\sigma_{a_m}\in {\cal C}_n$ and
$t_m,t_0\in \Delta$ such that $\lim_m t_m =t_0$ and $\lim_m \sigma_{a_m} (t_m)=\la_0$.
Moreover, $\sigma_{a_m}(\Delta)\subset Per_{k_m}(0)$ for some increasing sequence of integers $(k_m)_m$.
\item[2)]Let $w\in \Delta$. If $\la_0$ belongs to the bifurcation locus of $Per_n(w)$ then there exists $a_m\in ({\cal C}_n)_0$ and
$\sigma_{a_m}\in {\cal C}_n$ such that
$\lim_m \sigma_{ a_m}(w)=\la_0$. Moreover, $\sigma_{a_m}(\Delta)\subset Per_{k_m}(0)$ for some increasing sequence of integers $(k_m)_m$.
\end{itemize}
\end{lem}

\proof
 Let us consider the first assertion. Since $\la_0\in U_n$ there exists $t_0\in \Delta$
such that $p_n(\la_0,t_0)=0$. As $\la_0\in \textrm{Supp}(\Bc)$, it follows from corollary \ref{AccCent} that $\la_0=\lim_m\la_{k_m}$
where $\la_{k_m}\in U_n\cap Per_{k_m}(0)$ for some increasing sequence of integers $(k_m)_m$.
By continuity of $p_n$ there exists a sequence $(t_m)_m$ such that $\lim_m t_m=t_0$ and $p_n(\la_{k_m},t_m)=0$.
Now, since $U_n\cap Per_{k_m}(0) \subset H_{n,k_m}$, the point $\la_{k_m}$ belongs to some component $H_{n,k_m}^{j_m}$ and
we must have 
$\phi_{n,k_m}^{j_m}(\la_{k_m})=(t_m,0)$. According to the above definition, this exactly means that $\la_{k_m}=\sigma_{a_m}(t_m)$ where 
$a_m:=(\phi_{n,k_m}^{j_m})^{-1}(0,0)$. By construction $a_m\in ({\cal C}_n)_0$, $\sigma_{a_m}\in {\cal C}_n$ and $\sigma_{a_m}(\Delta)\subset Per_{k_m}(0)$.\\
 We now consider the second assertion. Since $\la_0 \in Per_n(w)$ we have
$p_n(\la_0,w)=0$.
As $\la_0$ is in the bifurcation locus of $Per_n(w)$ it follows from corollary \ref{AccCent} that $\la_0=\lim_n\la_{k_m}$
where $\la_{k_m}\in Per_n(w)\cap Per_{k_m}(0)$ for some increasing sequence of integers $(k_m)_m$.
Now, since $Per_n(w)\cap Per_{k_m}(0) \subset H_{n,k_m}$, the point $\la_{k_m}$ belongs to some component $H_{n,k_m}^{j_m}$ and 
$\phi_{n,k_m}^{j_m}(\la_{k_m})=(w,0)$. This exactly means that $\la_{k_m}=\sigma_{a_m}(w)$ where 
$a_m:=(\phi_{n,k_m}^{j_m})^{-1}(0,0)$. By construction $a_m\in ({\cal C}_n)_0$ and $\sigma_{a_m}(\Delta)\subset Per_{k_m}(0)$.\qed\\

We shall need the following fact for applying the extension theorem \ref{thext} to ${\cal C}_n$.
It will be obtained from lemma \ref{lemcent} after observing that any curve $Per_n(w)$ must have a non empty bifurcation locus.

\begin{cor}\label{Bif/cent}
Let $w\in \Delta$ and $Per_n^0(w)$ be a component of $Per_n(w)$. Then the set $Per_n^0(w)\cap ({\cal C}_n)_w$ is infinite.
\end{cor}

\proof

According to the second assertion of lemma \ref{lemcent} it suffices to show that the bifurcation locus of $Per_n^0(w)$ is not empty.
To this purpose we will observe that there exists $\la_0\in Per_n^0(w)$ which do have an unstable neutral cycle. 
We use the compactification of ${\cal M}_2={\bf C}^2$ discussed in subsection \ref{compa} and argue as follows.
Using Proposition \ref{Per/L}, one sees that if $\theta_0\in [0,2\pi[$ then the curves $Per_n^0(w)$ and ${Per}_1(e^{i\theta_0})$ do not intersect 
 on the line at infinity: 
${Per}_n^0(w)\cap{Per}_1(e^{i\theta_0})\cap{\cal L}=\emptyset$. These curves must 
therefore intersect in ${\bf C}^2$. Let us then pick a point $\la_0\in {Per}_n(w)\cap{Per}_1(e^{i\theta_0})\cap{\bf C}^2$
and a holomorphic family $(\varphi_u)_{u\in\Delta}$ in $Rat_2({\bf P}^1)$ such that
$\pi\circ\varphi(0)=\la_0$ and $\pi\circ\varphi(\Delta)\subset Per_n^0(w)$.
Since $\la_0\in{Per}_1(e^{i\theta_0})$, the map $\varphi_0$ has an unstable neutral cycle and therefore the bifurcation locus of the family
$(\varphi_u)_{u\in\Delta}$ is not empty.\qed

\subsection{Holomorphic motion and lamination 
in the bifurcation locus}

We are now in order to state and prove the main result of this section.

\begin{theo}\label{thMouvBif}
Let $U_n \subset {\cal M}_2 ={\bf C}^2$ be the open subset of parameters which 
do have an attracting cycle of period $n$.
Let $\Bif_n$ be the bifurcation locus in $U_n$ and $\Bc{\arrowvert_{U_n}}$ be 
the associated bifurcation current.
Let $\Bif_{n}^{\;c}$ be the bifurcation locus in the central curve $Per_n(0)$
and $\mu_{n}^c$ be the associated bifurcation measure.
Then, there exists a map
\begin{eqnarray*}
\sigma:\Bif_{n}^{\;c}\times\Delta & \longrightarrow &\Bif_n \\
(\la,t)&\longmapsto&\sigma_\la(t)\\
\end{eqnarray*}
such that:

\begin{itemize}
\item[1)] $\sigma\left(\Bif_{n}^{\;c}\times\Delta\right)=\Bif_n$
\item[2)] $\sigma$\; \textrm{\it is continuous,} $\sigma(\la,\cdot)\; \textrm{\it is one-to-one and holomorphic for each }\la \in \Bif_{n}^{\;c}$  
\item[3)] $p_n(\sigma_{\lambda}(t),t)=0;\;\forall \la\in \Bif_{n}^{\;c},\forall t\in \Delta$
\item[4)] \textrm{\it the discs} 
$\left(\sigma_\la (\Delta)\right)_{\la\in\Bif_{n}^{\;c}}$ 
\textrm{\it are mutually disjoint}.
\end{itemize}

Moreover the bifurcation current in $U_n$ is given by
$$\Bc{\arrowvert_{U_n}} = \int_{\Bif_{n}^{\;c}}[\sigma_\la(\Delta)]\;\mu_{n}^c$$

and, in particular, $ \Bif_n$  is a lamination with $ \mu_{n}^c$ as transverse measure. 
\end{theo}

Let us recall that the map $ \sigma:\Bif_{n}^{\;c}\times\Delta \longrightarrow \Bif_n$ is what we called a 
$p_n$-guided holomorphic motion (see section \ref{HolMo}).\\

\proof
We start with  the $p_n$-guided holomorphic motion of the 
centers ${\cal C}_n$ which was constructed in the former subsection (see Definition \ref{Mouvcenter}).\\

\underline{Ste}p\underline{ 1}. {\it Extension of the holomorphic motion of centers: there exists
a continuous $p_n$-guided holomorphic motion $\widehat{\cal C}_n$ such that
${\cal C}_n\subset\widehat{\cal C}_n$ and $(\widehat{\cal C}_n
)_0
=\overline{\left({\cal C}_n\right)_0}$.}\\

According to the corollary \ref{Bif/cent}, every component of the curve 
$\{p_n(\cdot,t)=0\}$ contains at least three points of 
$\left({\cal C}_n\right)_t$. Moreover, the degree of the polynomials $p_n(\cdot,t)$ does not depend on
$t\in\Delta$ (see theorem \ref{MilParam}). 
This will allow us to use theorem \ref{thext}. Let us choose an exhaustion
$\left(K_l\right)_l$ of
${\bf C}^2$ by compact subsets and consider the following subfamilies of ${\cal C}_n$:
$${\cal F}_l:=\{\sigma\in {\cal C}_n/\sigma(0)\in K_l\}.$$
We may apply theorem \ref{thext} to any pair $({\cal F}_l,{\cal C}_n)$. 
This gives an increasing sequence of $p_n$-guided holomorphic motions 
$(\widehat{\cal F}_l)_l$ and it then suffices to set
$$\widehat{\cal C}_n:=\cup_{l}\widehat{\cal F}_l.$$

\underline{Ste}p\underline{ 2}. {\it Restriction of $\widehat{\cal C}_n$ to $\Bif_{n}^{\; c}$.}\\

According to the second assertion of lemma \ref{lemcent} we have 
$\Bif_{n}^{\;c}\subset \overline{\left({\cal C}_n\right)_0}$ and, by the above first step this
also means that $\Bif_{n}^{\;c}\subset(\widehat{\cal C}_n)_0$. All we want to show here is that the 
collection of discs of $\widehat{\cal C}_n$ which go through $\Bif_{n}^{\;c}$ coincides with
$\Bif_n$:
$$\Bif_n=\bigcup_{\la\in \Bif_{n}^{\;c}}\sigma_\la(\Delta).$$

Let us pick $\la_1:=\sigma_{\la_0}(t_0)$ where $\la_0\in\Bif_{n}^{\;c}$ and $t_0\in \Delta$.
We will show that $\la_1\in \Bif_n$. By the second assertion of lemma \ref{lemcent}, 
$\la_0=\lim_na_m$ with $a_m\in ({\cal C}_n)_0$. Then, using the continuity of $\widehat{\cal C}_n$,
one gets $\la_1=\sigma_{\la_0}(t_0)=\lim_n\sigma_{a_m}(t_0)$. Since $\sigma_{a_m}(\Delta)\subset Per_{k_m}(0)$
for some increasing sequence of integers $(k_m)_m$, this implies that $\la_1\in \Bif_n$.\\
Let us pick $\la_0\in  \Bif_n$. By the first assertion of lemma \ref{lemcent} we have 
$\la_0=\lim_m \sigma_{a_m}(t_m)$ where $a_m\in ({\cal C}_n)_0$ and $\lim_n t_m=t_0\in \Delta$.
Let us consider ${\cal F}:=\{\sigma_{a_m};\; m\in\NN\}$, this is a subfamily of ${\cal C}_n$
and the pair $({\cal F}, {\cal C}_n)$ satisfies the assumptions of lemma \ref{lemnorma}.
Thus, after taking a subsequence, we may assume that $\sigma_{a_m}\to\sigma$ for some 
$p_n$-guided holomorphic disc $\sigma$. In particular, we have
$\lim_m a_m=\lim_m\sigma_{a_m}(0)=\sigma(0)=:a_0$ and, since   
$\sigma_{a_m}(\Delta)\subset Per_{k_m}(0)$
for some increasing sequence of integers $(k_m)_m$, this implies that $a_0 \in \Bif_{n}^{\;c}$.
We may now use the continuity of $\widehat{\cal C}_n$ and get:
$\sigma_{a_0}=\lim_m \sigma_{a_m}$ so that finally $\la_0=\lim_m \sigma_{a_m}(t_m)=\sigma_{a_0}(t_0) \in 
\cup_{\la\in \Bif_{n}^{\;c}}\sigma_\la(\Delta).$\\

\underline{Ste}p\underline{ 3}. {\it Laminarity of} $\Bc\vert_{U_n}$.\\

According to the approximation formula given by theorem \ref{appT} applied on the central curve $Per_n(0)$ we have:
\begin{eqnarray}\label{z}
\mu_{n}^c=\lim_m 2^{-m}\sum_{\la\in Per_n(0)\cap Per_m(0)}\delta_{\sigma_{\la}(0)}.
\end{eqnarray}
Let us set $T:=\int[\sigma_{\la}(\Delta)]\;\mu_{n}^c$. We have to check that $T=\Bc{\arrowvert_{U_n}}$.
Let $\phi$ be a $(1,1)$-test form in $U_n$. As the holomorphic motion $\sigma$ is continuous, the function
$\la\mapsto \langle[\sigma_{\la}(\Delta)],\phi\rangle$ is continuous as well. Then, using \ref{z} one gets
\begin{eqnarray}\label{zz}
\langle T,\phi\rangle =\lim_m 2^{-m}\sum_{\la\in Per_n(0)\cap Per_m(0)}\langle[\sigma_{\la}(\Delta)],\phi\rangle =
\lim_m2^{-m}\langle [Per_m(0)],\phi\rangle
\end{eqnarray}
where the last equality uses the fact that, according to proposition \ref{multi}, 
 the curves $Per_m(0)$ have no multiplicity in $U_n$. 
Now the conclusion follows by using \ref{zz} and the approximation formula of theorem \ref{appT} in $U_n$.\qed

\section{Motion from the quadratic       polynomial family}\label{U1}

Two fundamental facts in the Ma\~n\'e-Sad-Sullivan theory are the density of the set ${\cal S}$ of stable maps 
(the complement of $\Bif$) and the concept of hyperbolic
components: two elements lying in the same connected component of $\cal S$ are either both hyperbolic or both non-hyperbolic. This reduces 
Fatou's conjecture on the density of hyperbolic rational maps to the problem of existence of non-hyperbolic components. This question is still open, 
even for 
the family of quadratic polynomials. We refer to \cite{MacMu}, \cite{MacMu2} for more details and further results.\\
In this section, combining our previous results with Slodkowski's theorem, we obtain sharper statements for the region $U_1=
\{\la\in{\bf C}^2/ f_{\la} \;\textrm{has an attracting fixed point}\}$.
In particular, we prove that non-hyperbolic components exist in $U_1$ if and only if such components exist within the quadratic 
polynomial family $Per_1(0)$. Our result may be stated as follows
(we recall that the bifurcation locus in $Per_1(0)$ is the boundary of the Mandelbrot set $M_2$)

\begin{theo}\label{farf}
There exists a holomorphic motion
${\tilde\sigma}:Per_1(0)\times\Delta \longrightarrow U_1$ which is $p_1$-guided, continuous, onto and such that
${\tilde\sigma}\left(bM_2\times\Delta\right)=\Bif_1$. All stable components 
in $U_1$ are of the form 
${\tilde\sigma}\left(\omega\times\Delta\right)$ for some component $\omega$ in 
$Per_1(0)$.
Moreover, the map $\la\mapsto{\tilde\sigma}(\la,t)$ is a
quasi-conformal homeomorphism for each $t$ and
$\tilde\sigma$ is one-to-one on $\left( Per_1(0)\setminus{\overline\car}\right)\times \Delta$ where $\car$ is the 
 main cardioid.  
\end{theo}

The key point here is that $U_1=\cup_{t\in \Delta}Per_1(t)$ where each $Per_1(t)$ is a straight line in ${\bf C}^2$ (see proposition \ref{Per/L}). 
This allows to see the holomorphic motion $\sigma$ constructed in theorem \ref{thMouvBif} as given by holomorphic motion in ${\bf P}^1$ and use Slodkovsky 
theorem to extend it to all $Per_1(0)$. Our reference for holomorphic motions in ${\bf P}^1$ is the book \cite{Hub}.\\

\proof

We first show the existence of a $p_1$-guided holomorphic motion $\tilde{\sigma}:Per_1(0)\times\Delta  \longrightarrow U_1 $ which is continuous, onto and 
satisfies $\tilde{\sigma}\left(bM_2\times \Delta\right)=\Bif_1.$\\
According to theorem \ref{thMouvBif} there exists a $p_1$-guided holomorphic motion  
\begin{eqnarray*}
\sigma:\overline{\left({\cal C}_1\right)_0}\times\Delta & \longrightarrow & U_1 \\
(\la,t)&\longmapsto&\sigma_\la(t)\\
\end{eqnarray*}

which is continuous and satisfies $\tilde{\sigma}\left(bM_2\times\Delta\right)=\Bif_1$; we recall that $bM_2\subset \overline{\left({\cal C}_1\right)_0}$.
We thus have to exhibit an extension ${\tilde\sigma}:Per_1(0)\times\Delta \longrightarrow U_1$ of $\sigma$ with the property that
${\tilde\sigma}(\cdot,t):Per_1(0) \longrightarrow Per_1(t)$ is onto for all $t\in\Delta$.
\\ 

Let us recall that, by definition, we have
\begin{eqnarray}\label{inject}
\sigma_{\la}(\Delta)\cap\sigma_{\la'}(\Delta)=\emptyset\;\;\;\textrm{if}\;\;\la\ne\la'
\end{eqnarray}

and that by construction (see proposition \ref{Mouvcenter}):
\begin{eqnarray}\label{motocentri}
\forall t\in \Delta,\;\forall \la\in Per_m(0)\cap Per_1(t),\; \exists \la_0\in Per_m(0)\cap Per_1(0)\;\textrm{such that}\; \la=\sigma_{\la_0}(t).
\end{eqnarray}

Let us write $\sigma_{\la}(t)=:\left(\alpha_{\la}(t),\beta_{\la}(t)\right)$. As $\sigma$ is $p_1$-guided we have $\sigma_{\la}(t)\in Per_1(t)$ and thus, according to proposition \ref{Per/L}:
\begin{eqnarray}\label{equa}
\alpha_{\la}(t)= \frac{1}{1+t^2}\left(t\beta_{\la}(t)+t^3+2\right),\;\;\forall t\in \Delta.
\end{eqnarray}

Setting $\beta_{\infty}(t)=\infty$ for all $t\in\Delta$, we may consider the map $\beta$ as a holomorphic motion
in ${\bf P}^1$ that is $\beta:\overline{\left({\cal C}_1\right)_0}\cup\{\infty\}\times\Delta\longrightarrow {\bf P}^1$. In particular the injectivity
in $\la$ for fixed $t$ follows from \ref{inject} and \ref{equa}.\\
By Slodkowski's theorem, $\beta$ extends to some holomorphic motion of ${\bf P}^1$, that is, to some
 map $\tilde{\beta}:{\bf P}^1\times\Delta\longrightarrow {\bf P}^1$
which, in particular, is continuous and such that $\la\mapsto\tilde{\beta}_{\la}(t)$ is a (quasi-conformal) homeomorphism of ${\bf P}^1$ for
every $t\in \Delta$. By construction we have $\tilde{\beta}_{\la}(t)\ne\infty$ for $\la\ne\infty$ and, considering $\tilde{\beta}_{\la}(t)$ as a complex number
for $\la\in Per_1(0)$, it suffices to set:
$$\tilde{\sigma}(\la,t)=\left(t\tilde{\beta}_{\la}(t)+t^3+2,\tilde{\beta}_{\la}(t)\right).$$

We now want to show that all stable components in $U_1$ are of the form
${\tilde \sigma}\left(\omega\times\Delta\right)$ for some stable component
$\omega$ in $Per_1(0)$. We treat the main cardioid $\car$ before, because the 
argument we need is specific. Let us show that 
\begin{eqnarray}\label{cardio}
{\tilde \sigma}\left({\overline\Car}\times\Delta\right)=\overline{A_1}\;\;
\textrm{and}\;\;
{\tilde \sigma}^{-1}\left(\overline{A_1}\right)=
{\overline\Car}\times\Delta
\end{eqnarray}
where $A_1$ is the open subset of $U_1$ which consists of parameters possessing 
 two distinct attracting fixed points ($\car=A_1\cap
Per_1(0)$) and $\overline{A_1}$ its closure in $U_1$.
 
Let us first show that
${\tilde \sigma}\left(\Car\times\Delta\right)\subset\overline{A_1}$.
If this would
not be the case, we would find $\la_1\in\Car$ and $t_1\in\Delta$ such that
${\tilde\sigma}_{\la_1}(t_1)\notin\overline{A_1}$. Then, 
$\la_1={\tilde\sigma}_{\la_1}(0)$ and ${\tilde\sigma}_{\la_1}(t_1)$
  do not have the same number of attracting fixed
points and therefore, a bifurcation occurs within the holomorphic family 
$\{{\tilde\sigma}_{\la_1}(t);\;t\in\Delta\}$. Thus, by corollary 
\ref{AccCent}, there
exists $m\in\NN^*$ and $t_m\in\Delta$ such that 
  ${\tilde\sigma}_{\la_1}(t_m)\in Per_m(0)$. As $\tilde\sigma$ extends $\sigma$,
we then get from \ref{motocentri} that 
${\tilde\sigma}_{\la_1}(t_m)={\tilde\sigma}_{\la_0}(t_m)$ for some $\la_0\in
Per_1(0)\cap Per_m(0)$. 
As the map $\la\mapsto {\tilde\sigma}_{\la}(t_m)$ is one-to-one,
we have $\la_1=\la_0$. This is a contradiction since $\la_0$ does not belong to $\Car$.\\
As $\tilde\sigma$ is continuous we also have 
${\tilde \sigma}\left({\overline\Car}\times\Delta\right)\subset\overline{A_1}$ and it thus remains to show that
${\tilde \sigma}^{-1}\left(\overline{A_1}\right)\subset
{\overline\Car}\times\Delta$. We will essentially argue as before. Assume to the contrary that 
${\tilde\sigma}_{\la_1}(t_1)\in\overline{A_1}$ for some $\la_1\notin{\overline\Car}$ and $t_1\in\Delta$. 
Then ${\tilde\sigma}_{\la_1}(t_1)$ has a non-repulsiv fixed point which differs from $\infty$ but this is not the case for $\la_1$. 
A bifurcation must therefore occur within the holomorphic family 
$\{{\tilde\sigma}_{\la_1}(t);\;t\in\Delta\}$. As before,
there
exists $m\in\NN^*$ and $t_m\in\Delta$ such that 
  ${\tilde\sigma}_{\la_1}(t_m)\in Per_m(0)$ and, since $\tilde\sigma$ extends $\sigma$,
we get from \ref{motocentri} that 
${\tilde\sigma}_{\la_1}(t_m)={\tilde\sigma}_{\la_0}(t_m)$ for some $\la_0\in
Per_1(0)\cap Per_m(0)$. 
As the map $\la\mapsto {\tilde\sigma}_{\la}(t_m)$ is one-to-one,
we have $\la_1=\la_0$.
  This is impossible.
Indeed, as $\la_0\in Per_m(0)$ we would have 
${\tilde\sigma}_{\la_1}(\Delta)={\tilde\sigma}_{\la_0}(\Delta)=
{\sigma}_{\la_0}(\Delta)\subset Per_m(0)$ and
 ${\tilde\sigma}_{\la_1}(t_1)$  would have three
distinct non-repulsiv cycles (two fixed points and a $m$-cycle).
This contradicts the Fatou-Shishikura inequality.\\

To end the proof we will show that 
${\tilde\sigma}:\left(Per_1(0)\setminus{\overline\Car}\right)\times\Delta
\longrightarrow U_1\setminus\overline{A_1}$ is a
 homeomorphism. Since
${\tilde\sigma}(bM_2)=\Bif_1$, this will imply that
${\tilde\sigma}\left(\omega\times\Delta\right)$ is a stable component
in $U_1$ for every stable component $\omega\subset Per_1(0)$ which differs from $\car$.\\
As the map ${\tilde\sigma}$ is onto, it follows from 
\ref{cardio} that
${\tilde\sigma}:\left(Per_1(0)\setminus{\overline\Car}
\right)\times\Delta
\longrightarrow U_1\setminus\overline{A_1}$ is onto as well.
Let us check that it is one-to-one and open.
 If ${\tilde\sigma}_{\la}(t)={\tilde\sigma}_{\la'}(t')
\in U_1\setminus{\overline A_1}$ then, as two lines $Per_1(t)$ and $Per_1(t')$
cannot intersect in $  U_1\setminus{\overline A_1}$ one gets $t=t'$, and
$\la=\la'$ follows 
from the injectivity of $\la\mapsto {\tilde\sigma}_{\la}(t)$.\\
To establish the openess, one may use the fact that 
$U_1\setminus{\overline A_1}$ is foliated by the straight lines $Per_1(t)$ and 
that  $\la\mapsto {\tilde\sigma}_{\la}(t)\in Per_1(t)$ is a $K$-quasi-conformal
homeomorphism for $K:=\frac{1+\vert t\vert}{1-\vert t\vert}$.\qed

\bibliographystyle{amsalpha}

\end{document}